\documentclass[a4paper,10pt]{article}

\usepackage{authblk}

\usepackage[utf8]{inputenc} 
\usepackage[T1]{fontenc}    
\usepackage{indentfirst}
\usepackage{lmodern}

\usepackage{mathtools}
\usepackage{verbatim}
\usepackage{url}            
\usepackage{booktabs}       
\usepackage{amsfonts}        
\usepackage{amsmath}
\usepackage{amssymb}
\usepackage{amsthm}
\usepackage{microtype}      
\usepackage[a4paper]{geometry}
\usepackage{enumitem}
\usepackage{xcolor}
\usepackage{mathrsfs}

\usepackage{geometry}
\usepackage{fancyhdr}
\usepackage{listings}

\usepackage[linktocpage=true]{hyperref}
\hypersetup{
    colorlinks=true,
    linkcolor=blue,
    citecolor=violet,
    urlcolor=blue
}

\usepackage{graphicx}
\usepackage[makeroom]{cancel}
\usepackage{tikz}

\usepackage[inkscapeformat=png]{svg}
\usepackage[normalem]{ulem}

\theoremstyle{plain}
\newtheorem{theorem}{Theorem}[section]

\theoremstyle{definition}

\theoremstyle{remark}

\makeatletter
\def\keywords{\xdef\@thefnmark{}\@footnotetext}
\makeatother

\title{Probability measure annihilating all finite-dimensional subspaces}
\author{Nizar El Idrissi and Hicham Zoubeir}

\newcommand{\Addresses}{{
  \bigskip
  \footnotesize
  \textbf{Nizar El Idrissi.}
  \par\nopagebreak 
  Laboratoire : Equations aux dérivées partielles, Algèbre et Géométrie spectrales.
  \par\nopagebreak 
  Département de mathématiques, faculté des sciences, université Ibn Tofail, 14000 Kénitra.
  \par\nopagebreak 
  \textit{Adresse e-mail} : \texttt{nizar.elidrissi@uit.ac.ma}

  \medskip
  \footnotesize
  \textbf{Hicham Zoubeir.}
  \par\nopagebreak 
  Laboratoire : Equations aux dérivées partielles, Algèbre et Géométrie spectrales.
  \par\nopagebreak 
  Département de mathématiques, faculté des sciences, université Ibn Tofail, 14000 Kénitra.
  \par\nopagebreak 
  \textit{Adresse e-mail} : \texttt{hzoubeir2014@gmail.com}
}}

\begin{document}
\newpage
\maketitle

\begin{abstract}
We propose in this short note a prime numbers-based method for constructing probability measures on infinite-dimensional Banach spaces annihilating all finite-dimensional subspaces, supplementing the methods of construction of Gaussian measures and infinite-product-type probability measures. This new method confirms that probability measures with this property are generic amongst probability measures that are not supported on finite-dimensional subspaces. In the process, we show the existence of an uncountable measurable family of independent vectors having the cardinality of the continuum in any infinite-dimensional Banach space.
\end{abstract}

\keywords{2020 \emph{Mathematics Subject Classification.} 46B09; 28A05; 03E10; 15A03; 11N05; 13A05}
\keywords{\emph{Key words and phrases.} Banach space; probability measure; cardinality of the continuum; linear independence; primes distribution; unique factorization property}

\tableofcontents

\section{Introduction}

The structure of probability measures on infinite-dimensional Banach spaces lies at the crossroads of functional analysis and probability theory. A fundamental dichotomy distinguishes measures supported on finite-dimensional subspaces from genuinely infinite-dimensional measures assigning measure zero to every finite-dimensional vector subspace. Such measures play an important role in the study of linear independence of random vectors, genericity phenomena, and probabilistic constructions in Banach and Hilbert spaces (see e.g. \cite{VakhaniaTarieladzeChobanyan2012} and \cite{Kukush2020}). \\ \\
In our previous work \cite{ElIdrissiZoubeir2026} (following results in \cite{ChristensenHasannasab2019} and \cite{ElIdrissiKabbaj2023}), we showed that if a sequence of Hilbert-space-valued random vectors has finite-dimensional marginals absolutely continuous with respect to a probability measure that annihilates all finite-dimensional subspaces, then the sequence is almost surely linearly independent. That result highlights the structural importance of such measures. However, while classical constructions (Gaussian measures, infinite product measures) provide examples, the existence of probability measures annihilating all finite-dimensional subspaces deserves a direct treatment. \\ \\
The first objective of this paper is to provide a new and explicit construction of such measures on arbitrary infinite-dimensional Banach spaces. Our method is elementary but arithmetically structured: it uses the unique factorization of integers and asymptotic properties of prime numbers (namely, the Baker–Harman–Pintz theorem \cite{BakerHarman1996,BakerHarmanPintz2001}) to build a continuum-sized family of infinite subsets of $\mathbb{N}$ with pairwise finite intersections \cite{Buddenhagen1971} and suitable measurability properties. This number-theoretic encoding allows us to construct an injective measurable map $x : (0,1) \to E$ such that $(x_t)_{t \in (0,1)}$ is linearly independent, whose pushforward of the Lebesgue measure yields a Borel probability measure on $E$ that annihilates every finite-dimensional subspace. Therefore, we provide a new proof for:

\begin{theorem}
Let $E$ be an infinite-dimensional Banach space over $\mathbb{R}$ or $\mathbb{C}$. Then there exists a Borel probability measure $P$ on $E$ annihilating all finite-dimensional vector subspaces.
\end{theorem}

Beyond existence, the construction reveals a deeper structural fact. It produces an uncountable measurable family $(x_t)_{t \in (0,1)}$ of linearly independent vectors in $E$,  indexed by the continuum. While it is classical that the Hamel dimension of an infinite-dimensional separable Banach space is the continuum (see \cite{Lacey1973}), our result strengthens this perspective by ensuring that the family is obtained through a measurable parametrization compatible with the Borel structure of $E$.
More precisely, we prove:

\begin{theorem}
There is an uncountable measurable family of independent vectors having the cardinality of the continuum in any infinite-dimensional Banach space. 
\end{theorem}

The measurability aspects of the construction require some care and connect naturally with tools from descriptive set theory (see \cite{Kechris1995}). In particular, we establish the measurability of the range of the parametrizing map and show that its closure differs from the range by at most a countable set. \\ \\
Conceptually, the paper shows that probability measures annihilating finite-dimensional subspaces are not exceptional: they arise from flexible constructions available in every infinite-dimensional Banach space. The method of construction discussed in this article complements classical Gaussian and product constructions and provides a new viewpoint linking prime number theory, linear independence, and probability theory in infinite-dimensional spaces.

\section{Proof of both theorems}

Let $E$ be an infinite-dimensional Banach space over $\mathbb{R}$ or $\mathbb{C}$ and denote by $\mathcal{P}$ the set of prime numbers. The following proof is motivated by the fact that there is no atomless probability measure on a countable set. Since $(0,1) \subset \mathbb{R}$ admits such a probability measure $\lambda$ (the Lebesgue measure for instance), we will define the probability measure $P$ on $E$ as a pushforward of $\lambda$ by an injective and measurable map $x : (0,1) \to E$. \\ \\
The first and second parts of the proof are inspired by the note \cite{Lacey1973}. However, instead of an arbitrary choice of $(N_t)_{(0,1)}$ such that
\begin{equation}
\label{NtConditions}
\forall t \in (0,1) : N_t \text{ is infinite and } \forall 0 < t \neq t' < 1 : N_t \cap N_{t'} \text{ is finite},
\end{equation}
we provide a specific construction of $(N_t)_{(0,1)}$ satisfying not only (\ref{NtConditions}) but also some measurability conditions. The third part of the proof is the smallest one and concludes the proof. \\

\begin{itemize}

\item \textbf{First part.} \\
For all $t \in (0,1)$ and $j \in \mathbb{N}^*$ such that $2 \leq t^{-j}$, consider
\[ e_j(t) := \max \{ p \in \mathcal{P} : p \leq t^{-j} \} \]
and if $t^{-j} < 2$
\[ e_j(t) := 2 \]
Further, define
\[ f_j(t) := \Pi_{i=1}^j e_i(t). \]
Clearly, $(e_j(t))$ and $(f_j(t))$ are nondecreasing in $j$ and nonincreasing in $t$. Moreover, $e_j(t)$ and $f_j(t)$ go to $+\infty$ as $t^{-j}$ goes to $+\infty$ since there are infinitely many primes. Furthermore, $(f_j(t))$ is increasing in $j$ since all prime numbers are greater or equal than 2. For all $t \in (0,1)$, consider
\[ N_t := \{ f_j(t) : j \in \mathbb{N}^* \} \subseteq \mathbb{N}^*. \]
Clearly $N_t$ is infinite for all $t \in (0,1)$. Let's show that $N_t \cap N_{t'}$ is finite for all $t ` > t$. Fix $0 < t < t' < 1$. Let $C,A,\alpha$ be three positive constants, with $\alpha \in (0,1)$ such that
\[ t^{-A} \geq t'^{-A} \geq 2, \quad \quad (*) \]
\[ \forall x \in [t'^{-A},\infty( \; : \quad \quad x - C x^\alpha \le p(x), \quad \quad (**) \]
\[ \forall x \in [A,\infty( \; : \quad \quad C t^{(1-\alpha)x} < \frac{1}{2}, \quad \quad (***) \]
and
\[ \forall x \in [A,\infty( \; : \quad \quad (\frac{t}{t'})^x < \frac{1}{2}, \quad \quad (****) \]
where the inequality (**) is provided for instance by the Baker–Harman–Pintz theorem \cite{BakerHarman1996,BakerHarmanPintz2001}. Here $p(x) := \max\{p \in \mathcal{P} : p \leq x\}$ for all $x \geq 2$. Suppose by contradiction that $m \in N_t \cap N_{t'}$ with $m \geq t^{-A^2}$. Then there exists two natural numbers $u,v$ such that $m = f_u(t) = f_v(t')$, implying that $u=v$ and $e_i(t) = e_i(t')$ for $i$ in $[\![1,u]\!]$, since $f_u(t)$ and $f_v(t')$ are a product of $u$ and $v$ prime numbers respectively (this is the unique factorization property in $\mathbb{N}^*$). Since $t^{-A^2} \le m = f_u(t) \le (e_u(t))^u \le t^{-u^2}$, we have $u > A$. Now, by plugging $x := t^{-u}$ in inequality (**) and taking in account (*), we have
\[ t^{-u} - C t^{-\alpha u} \leq p(t^{-u}) = e_u(t) = e_u(t') \]
and so, multiplying by $t^u$ and using the inequalities (***) and (****),
\[ 1 \le C t^{(1-\alpha)u} + t^u e_u(t') \leq C t^{(1-\alpha)u} + (\frac{t}{t'})^u t'^u e_u(t') \le C t^{(1-\alpha)u} + (\frac{t}{t'})^u < 1. \]
This is a contradiction, proving that $N_t \cap N_{t'}$ is finite.

\item \textbf{Second part.} \\
Using the Hahn-Banach theorem multiple times and a variant of the Gram-Schmidt process, we can easily prove that there exists a sequence $(x_j,x_j^*)_{j \in \mathbb{N}^*}$ such that $(x_j,x_j^*) \in E \times E^*$ for all $j \in \mathbb{N}^*$, such that $x_u$ has norm 1 for all $u \in \mathbb{N}^*$, $x_v^*(x_u) = 0$ for all $u \neq v \in \mathbb{N}^*$ and $x_u^*(x_u) \neq 0$ for all $u \in \mathbb{N}^*$. For all $t \in (0,1)$, define 
\[ x_t = \sum_{ j\in N_t} \frac{x_j}{2^j}. \]
The construction of $(N_t)_{(0,1)}$ satisfying (\ref{NtConditions}) was for instance given in \cite{Buddenhagen1971}. However, the measurability of $x : (0,1) \to E$ and $Range(x)$ was not established. The construction of the $(N_t)_{(0,1)}$ that was provided in part 1 ensures these measurability properties, as the discussions here and in the first appendix show. First, observe as in \cite{Lacey1973} that $(x_t)_{(0,1)}$ is a linearly independent family in $E$ having the cardinality of the continuum and that $x : (0,1) \to x_t \in E$ is injective. Denote 
\[ S := \{ p^{-\frac{1}{i}} : p \in \mathcal{P}, i \in \mathbb{N}^* \} \]
For any $p \in \mathcal{P}$ and $i \in \mathbb{N}^*$, we have
\[ (e_i)^{-1}(\{p\}) = ((p_\text{succ})^{-\frac{1}{i}}, p^{-\frac{1}{i}}] \]
where we set $2_{succ} := \infty$ to account for the condition $e_i(t) = 2$ if $t^{-i} < 2$. In particular, for any $t^* \in (0,1) \setminus S$, we have
\[ \forall i \in \mathbb{N}^* : (e_i)^{-1}(\{e_i(t^*)\}) = ((e_i(t^*))_\text{succ}^{-\frac{1}{i}} , e_i(t^*)^{-\frac{1}{i}} ) \]
and so for all $i \in \mathbb{N}^*$, $e_i$ is continuous at all the points in $(0,1) \setminus S$. In general, for all $i \in \mathbb{N}^*$, $e_i$ is measurable as a function from $(0,1)$ to $\mathbb{N}$, and the same goes for the $f_i$. The fact that $x : (0,1) \to E$ is measurable is a straightforward consequence of the fact that the $f_j(t)$ are measurable in $t \in (0,1)$ for all $j \in \mathbb{N}$, plus the obvious formula
\[ x_t = \sum_{ u \in \mathbb{N}} \frac{x_{f_u(t)}}{2^{f_u(t)}}. \]

\item \textbf{Third part.} \\
For any measurable subset $M$ of $E$, set
\[ P(M) := \lambda(x^{-1}(M)). \]
Clearly, $P$ is a probability measure on $E$ (the pushforward of $\lambda$ by $x$). Moreover, if $M$ is a finite-dimensional subspace of $E$, the set $M \cap Range(x)$ is finite and so $P(M)=\lambda(x^{-1}(M \cap  Range(x))) = 0$ since $x$ is injective and $\lambda$ is atomless. Therefore, the proof is finished \qed.

\end{itemize}

\appendix

\section{Measurability of Range(x)}

To show the measurability of Range(x), we could simply invoke corollary 15.2 p. 89 of \cite{Kechris1995}. However, since this corollary was not known to us during the writing of our paper, we take the time in this appendix to prove this result on our own. Actually, we show the stronger property that Range(x) is the set difference $F \setminus C$, where $F$ is a closed set and $C$ is a countable set. Combining the measurability of $x : (0,1) \to E$ with the measurability of its range (which is really just its consequence), this implies that there is an uncountable measurable family of independent vectors having the cardinality of the continuum in any infinite-dimensional Banach space. As said in the introduction, this improves on the result of \cite{Lacey1973} and paves the way for further strengthenings. \\ \\
We prove that $Range(x)$ is a measurable subset of $E$. \\
Certainly, $\overline{Range(x)}$ is a measurable subset of $E$ since it is closed, so it suffices to show that $\overline{Range(x)} \setminus Range(x)$ is countable. \\
Consider the following cases:
\begin{itemize}
\item Let $t^* \in (0,1) \setminus S$ and consider
\[ X_{t^*} := \{l \in E : (\forall m \in \mathbb{N}^*)(\forall \eta > 0)(\exists t \in (t^*-\eta, t^* + \eta)) : \lVert x_t - l \rVert \leq \frac1m\} \]
From the continuity of $x : (0,1) \to E$ at $t^* \in (0,1) \setminus S$, we have that $X_{t^*} = \{x_{t^*} \}$.

\item Consider
\[ X_0 := \{l \in E : (\forall m \in \mathbb{N}^*)(\forall \eta > 0)(\exists t \in (0, \eta)) : \lVert x_t - l \rVert \leq \frac1m\} \]
Since 
\[ \lVert x_t \rVert = \lVert \sum_{j \in \mathbb{N}^*} \frac{x_{f_j(t)}}{2^{f_j(t)}} \rVert  \leq \sum_{u \geq f_1(t)} \frac{1}{2^u} = \frac{1}{2^{f_1(t)-1}} \to 0 \]
as $t \to 0^+$, this shows that $X_0 = \{0\}$.

\item Consider
\[ X_1 := \{l \in E : (\forall m \in \mathbb{N}^*)(\forall \eta > 0)(\exists t \in (1-\eta,1)) : \lVert x_t - l \rVert \leq \frac1m\} \]
For all $u \in \mathbb{N}^*$, the family of natural numbers $f_u(t)$ (in $t \in (0,1)$) converges to $2^u$ as $t \to 1^-$, so it is equal to $2^u$ for $t$ close enough to $1^-$. \\
Therefore
\[ x_t = \sum_{ u \in \mathbb{N}} \frac{x_{f_u(t)}}{2^{f_u(t)}} \to \sum_{u \in \mathbb{N}^*} \frac{x_{2^u}}{2^{2^u}} \]
as $t \to 1^-$ by interchanging the limit with the sum. \\
This shows that $X_1 = \{  \sum_{u \in \mathbb{N}^*} \frac{x_{2^u}}{2^{2^u}} \}$.

\item Let $t^* \in S$ and consider 
\[ X_{t^*}^- := \{l \in E : (\forall m \in \mathbb{N}^*)(\forall \eta > 0)(\exists t \in (t^*-\eta,t^*)) : \lVert x_t - l \rVert \leq \frac1m\} \]
From the left continuity of $x : (0,1) \to E$ at $t^* \in S$, we have that $X_{t^*} = \{x_{t^*} \}$.

\item Let $t^* = p^{-\frac{1}{i}} \in S$ and consider
\[ X_{t^*}^+ := \{l \in E : (\forall m \in \mathbb{N}^*)(\forall \eta > 0)(\exists t \in (t^*,t^*+\eta)) : \lVert x_t - l \rVert \leq \frac1m\} \]
For all $u \in \mathbb{N}^*$, the nonincreasing family of natural numbers $f_u(t)$ (in $t \in (0,1)$) converges to $f_u(t^*)^- < f_u(t^*)$ as $t \to (t^*)^+$, so it is equal to $f_u(t^*)^-$ for $t$ close enough to $(t^*)^+$. Therefore,
\[ x_t = \sum_{ u \in \mathbb{N}} \frac{x_{f_u(t)}}{2^{f_u(t)}} \to  \sum_{ u \in \mathbb{N}} \frac{x_{(f_u(t^*))^-}}{2^{(f_u(t^*))^-}} \]
as $t \to (t^*)^+$ by interchanging the limit with the sum (after noticing again that $f_u(t) \geq 2^u$, which allows for the limit-sum interchange). This shows that $X_{t^*}^+ = \{ \sum_{ u \in \mathbb{N}} \frac{x_{(f_u(t^*))^-}}{2^{(f_u(t^*))^-}} \}$.
\end{itemize}

Now, consider a sequence $x_{t_n}$ with $t_n \in (0,1)$ that converges to some $l \in E$. Since $[0,1]$ is compact, we can extract from $(t_n)$ a subsequence (that we still call $(t_n)$ by abuse of notation) such that $t_n \to t^* \in [0,1]$. So, we have
\[ \overline{Range(x)} \setminus Range(x) = X_0 \cup X_1 \cup \left( \bigcup_{t^* \in S}  (X_{(t^*)^+} \cup X_{(t^*)^-} \right). \]
and since $S$ is countable, this set is countable too. This proves by the previous discussion that $Range(x)$ is measurable.

\section*{Conflict of interest}
On behalf of all authors, the corresponding author states that there is no conflict of interest.  

\section*{Data availability statement}
No data are associated with this article.

\nocite{*}
\bibliographystyle{alpha}
\bibliography{references}

\Addresses

\end{document}